\documentclass{amsart} 
\usepackage{amssymb, amsmath}
\usepackage[dvips]{graphicx}
\usepackage{amsfonts}
\usepackage{latexsym}
\usepackage{color}

\newtheorem{theorem}{Theorem}	
\newtheorem{lemma}{Lemma}
\newtheorem{corollary}{Corollary}		
\newtheorem{proposition}{Proposition}		
\newtheorem{question}{Question}			
	
\newtheorem{definition}{Definition}

\title[Jacobian-squared function-germs]
{
Jacobian-squared function-germs   
}
\author{Takashi Nishimura
}
\address{
Takashi Nishimura: Research Institute of Environment 
and Information Sciences,  
Yokohama National University, 
Yokohama 240-8501, Japan}
\email{nishimura-takashi-yx@ynu.jp}
\begin{document}
\begin{abstract}
{
In this paper, 
it is shown that, for any equidimensional $C^\infty$ map-germ   
$f: (\mathbb{R}^n,0)\to (\mathbb{R}^n,0)$, the map-germ  
$F: (\mathbb{R}^n, 0) \to 
\mathbb{R}^n\times\mathbb{R}^{\ell}$ defined by 
$F(x)=\left(f(x), \mu_1(x){|Jf|^2(x)}, \cdots, 
\mu_\ell(x){|Jf|^2(x)}\right)$ is always a frontal; 
where $\mu_i$ is a $C^\infty$ function-germ and 
$|Jf|$ is the Jacobian-determinant of $f$.     
Moreover, it is also shown that when the multiplicity of $f$ is less than or 
equal to $3$, 
any frontal constructed from $f$ must be $\mathcal{A}$-equivalent to 
a frontal $F$ of the above form.        
} 
\end{abstract}
\maketitle  
\section{Introduction}\label{section 1} 
\noindent
Throughout this paper, let $n, \ell$ be positive integers; 
and all map-germs, vector fields and differential forms 
are of class $C^\infty$ unless otherwise stated.   
\par 
Let $f: (\mathbb{R}^n,0)\to (\mathbb{R}^n,0)$ 
 be a map-germ  and let $|Jf|$ denote the Jacobian-determinant of $f$.    
The square of $|Jf|$ is called 
the \textit{Jacobian-squared} function-germ of $f$.      
In the theory of singularities of mappings, 
it is well-known that the function-germ $|Jf|$ plays an essential role 
to investigate the behavior of $f$ 
(for instance, see \cite{arnoldguseinzadevarchenko, wall}).    
However, to the best of author's knowledge, so far 
there have been no literatures 
to emphasize the importance of the Jacobian-squared $|Jf|^2$.   
In this paper, it is explained that $|Jf|^2$ is very significant  
to construct non-trivial frontals from a given equidimensional 
map-germ $f: (\mathbb{R}^n,0)\to (\mathbb{R}^n,0)$.   
\par 
Let $F: (\mathbb{R}^n,0)\to (\mathbb{R}^{n+\ell}, 0)$ 
be a map-germ and let $T\mathbb{R}^{n+\ell}$ be the 
tangent bundle of $\mathbb{R}^{n+\ell}$.    
A map-germ $\Phi: (\mathbb{R}^n,0)\to T\mathbb{R}^{n+\ell}$ 
is called a \textit{vector field along $F$} if 
the equality $\pi\circ \Phi(x)=F(x)$ holds, where 
$\pi: T\mathbb{R}^{n+\ell}\to \mathbb{R}^{n+\ell}$ is the 
canonical projection.   
Namely, $\Phi$ has the form $\Phi(x)=(F(x), \phi(x))$ where 
$\phi(x)\in T_{F(x)}\mathbb{R}^{n+\ell}$.     
A map-germ $F: (\mathbb{R}^n,0)\to \mathbb{R}^{n+\ell}$ 
is called a \textit{frontal} if there exist vector fields 
$\Phi_1, \ldots, \Phi_\ell: (\mathbb{R}^n,0)\to T\mathbb{R}^{n+\ell}$ 
along $F$ such that the following three conditions are satisfied:   
\begin{enumerate}
\item[(1)] $\phi_i(x)\cdot tF(\xi)(x)=0$ for any $i$ $(1\le i\le \ell)$ 
and any $\xi\in \theta(n)$, 
where $\Phi_i(x)=\left(F(x), \phi_i(x)\right)$ and the dot in the center 
stands for the scalar product of two vectors in 
$T_{F(x)}\mathbb{R}^{n+\ell}$.   
\item[(2)] $\phi_i(0)\ne 0$ for any $i$ $(1\le i\le \ell)$.   
\item[(3)] $\phi_1(0), \ldots, \phi_\ell(0)$ are linearly independent.      
\end{enumerate} 
\noindent 
Here, $tF$ and $\theta(n)$ are the notations 
defined by J.~Mather in \cite{mather3} (For details, see \S \ref{section 2}).   
\par 
Though frontals have been already well-investigated 
(for instance see \cite{ishikawa}),      
it seems to be desired to obtain how to construct 
non-trivial frontals easily and systematically 
from a given equidimensional map-germ 
$f: (\mathbb{R}^n,0)\to (\mathbb{R}^n,0)$ which is not necessarily finite 
(for the definition of finite map-germ, see Subsection \ref{subsection 2.1}).   
\begin{theorem}\label{theorem 1}
Let $f: (\mathbb{R}^n,0)\to (\mathbb{R}^n,0)$ be a map-germ.    
For any $i$ $(1\le i\le \ell)$, let 
$\mu_i: (\mathbb{R}^n,0)\to \mathbb{R}$ 
be a function-germ.   
Then, the map-germ $F: (\mathbb{R}^n,0)\to \mathbb{R}^{n+\ell}$ 
defined by 
\[
F=\left(f, \mu_1|Jf|^2, \ldots, 
\mu_\ell|Jf|^2
\right)
\]   
is always a frontal.    
\end{theorem}
\begin{definition}[\cite{ishikawaaustralia, ishikawa}]
\label{ramification module}
{\rm 
For a given equidimensional map-germ 
$f=(f_1, \ldots, f_n): (\mathbb{R}^n,0)\to (\mathbb{R}^n,0)$, 
the \textit{ramification module} of $f$ (denoted by $\mathcal{R}_f$) 
is defined as the $f^*\left(\mathcal{E}_n\right)$-module consisting of 
all functions $\varphi$ such that $d\varphi$ is an element of 
the $\mathcal{E}_n$-module generated by $df_i$ $(i=1, \ldots, n)$, 
where $dh$ for a function-germ $h: (\mathbb{R}^n,0)\to \mathbb{R}$ 
stands for the exterior differential of $h$.  
}
\end{definition} 
\noindent 
Given an equidimensional map-germ 
$f: (\mathbb{R}^n,0)\to (\mathbb{R}^n,0)$, 
the \textit{pull-back} $f^*:\mathcal{E}_n\to \mathcal{E}_n$ is defined 
by $f^*(\eta)=\eta\circ f$.    The pull-back $f^*$ gives an  
$\mathcal{E}_n$-module structure via $f^*$.    
In this paper, an $\mathcal{E}_n$-module via $f^*$ is called and 
denoted by an $f^*\left(\mathcal{E}_n\right)$-module.     Hence, 
$f^*\left(\mathcal{E}_n\right)$ itself is 
an $f^*\left(\mathcal{E}_n\right)$-module.    
It is clear that $f^*\left(\mathcal{E}_n\right)$ itself is contained in $\mathcal{R}_f$.    
By using the notion of the ramification module $\mathcal{R}_f$, 
Theorem \ref{theorem 1} is equivalent to assert that 
\[
\langle |Jf|^2\rangle_{\mathcal{E}_n}+f^*\left(\mathcal{E}_n\right)
\;\subset\; \mathcal{R}_f, 
\]
where $\langle |Jf|^2\rangle_{\mathcal{E}_n}$ is the $\mathcal{E}_n$-module 
generated by 
$|Jf|^2$.        
\begin{definition}[\cite{ishikawaaustralia, ishikawa}]
\label{opening}
{\rm 
Let $f: (\mathbb{R}^n,0)\to (\mathbb{R}^n,0)$ be an equidimensional 
map-germ and 
let $\psi_1, \ldots, \psi_\ell$ be elements of $\mathcal{R}_f$.    
Then, the map-germ 
\[
(f, \psi_1, \ldots, \psi_\ell) : (\mathbb{R}^n,0)\to 
\mathbb{R}^n\times \mathbb{R}^\ell
\] 
is called an \textit{opening} of $f$.    
}
\end{definition}
\noindent 
By definition, for any equidimensional map-germ, 
its opening is always a frontal.         
Theorem \ref{theorem 1} guarantees that 
even for non-finite equidimensional map-germs, 
it is constructed automatically infinitely many non-trivial openings.   
This is one advantage of Theorem \ref{theorem 1}.      
In Subsections \ref{fold singularity}-\ref{opening of $4_k$}, 
by direct elementary calculations, 
it is easily shown that the normal forms of known frontals are actually 
$\mathcal{A}$-equivalent to the frontals constructed by 
Theorem \ref{theorem 1}.      
In other words, 
in Subsections \ref{fold singularity}-\ref{opening of $4_k$}, 
frontals constructed by using Jacobian-squared function-germs 
are already very near (with respect to $\mathcal{A}$-equivalence) 
to the normal forms of celebrated frontals, and 
therefore   
criteria for these noticeable frontals are not needed.     
Since the normal forms of these frontals seem to be not easy to memorize, 
this may be another advantage of Theorem \ref{theorem 1}.   
\begin{proposition}[\cite{ishikawa}]\label{proposition frontal}
For any frontal germ 
$F: (\mathbb{R}^n,0)\to (\mathbb{R}^{n+\ell}, 0)$, 
there exist germs of diffeomorphism 
$h: (\mathbb{R}^n,0)\to (\mathbb{R}^n,0)$ and 
$H: (\mathbb{R}^{n+\ell},0)\to (\mathbb{R}^{n+\ell},0)$, 
an equidimensional map-germ 
$f: (\mathbb{R}^n,0)\to (\mathbb{R}^n,0)$ and 
elements $\psi_1, \ldots, \psi_\ell$ of $\mathcal{R}_f$ 
such that the following equality holds: 
\[
H\circ F\circ h = (f, \psi_1, \ldots, \psi_\ell).   
\]
\end{proposition}
Based on Proposition \ref{proposition frontal}, it is natural to ask  
the converse of Theorem \ref{theorem 1}.     
However, it turns out that there exist counterexamples 
against the converse of 
Theorem \ref{theorem 1} (see Subsection \ref{$A_k$-front singularity}).    
Nevertheless, Subsections 
\ref{fold singularity}-\ref{opening of $4_k$} 
suggest that if $f$ satisfies 
$\dim_{\mathbb{R}}Q(f)\le 3$, then asking the converse of 
Theorem \ref{theorem 1} is reasonable (for the definition of $Q(f)$, see 
Section \ref{section 2}).     This question can be answered affirmatively 
as follows.   
\begin{theorem}\label{theorem 2}
Let $F: (\mathbb{R}^n,0)\to (\mathbb{R}^{n+\ell}, 0)$ be a frontal germ.    
Suppose that there exist germs of diffeomorphism 
$h: (\mathbb{R}^n,0)\to (\mathbb{R}^n,0)$ and 
$H: (\mathbb{R}^{n+\ell},0)\to (\mathbb{R}^{n+\ell},0)$, 
an equidimensional map-germ 
$f: (\mathbb{R}^n,0)\to (\mathbb{R}^n,0)$ with  
$\dim_{\mathbb{R}}Q(f)\le 3$ 
and 
elements $\psi_1, \ldots, \psi_\ell$ of $\mathcal{R}_f$ 
such that the following equality holds: 
\[
H\circ F\circ h = (f, \psi_1, \ldots, \psi_\ell).   
\]
Then, the following holds:   
\[
\left\langle |Jf|^2 \right\rangle_{\mathcal{E}_n}+f^*\left(\mathcal{E}_n\right)
=\mathcal{R}_f.   
\]
\end{theorem}
\begin{corollary}
\label{corollary 1}
Let $F: (\mathbb{R}^n,0)\to (\mathbb{R}^{n+\ell}, 0)$ be a frontal germ.    
Suppose that there exist germs of diffeomorphism 
$h: (\mathbb{R}^n,0)\to (\mathbb{R}^n,0)$ and 
$H: (\mathbb{R}^{n+\ell},0)\to (\mathbb{R}^{n+\ell},0)$, 
an equidimensional map-germ 
$f: (\mathbb{R}^n,0)\to (\mathbb{R}^n,0)$ with  
$\dim_{\mathbb{R}}Q(f)\le 3$ 
and 
elements $\psi_1, \ldots, \psi_\ell$ of $\mathcal{R}_f$ 
such that the following equality holds: 
\[
H\circ F\circ h = (f, \psi_1, \ldots, \psi_\ell).   
\]
Then,  there exist a germ of diffeomorphism 
$\widetilde{H}: 
(\mathbb{R}^{n+\ell},0)\to (\mathbb{R}^{n+\ell},0)$ 
and function-germs $\mu_i: (\mathbb{R}^n, 0)\to \mathbb{R}$ 
$(1\le i\le \ell)$ such that 
\[
\widetilde{H}\circ H\circ F\circ h =
 (f, \mu_1|Jf|^2, \ldots, \mu_\ell|Jf|^2).  
\]
\end{corollary}     
From Theorem \ref{theorem 2}, 
it is natural to ask the following:   
\begin{question}\label{question 1}
Let $f: (\mathbb{R}^n,0)\to (\mathbb{R}^n,0)$ be an equidimensional 
map-germ.   
Then, does there exist a finitely generated $\mathcal{E}_n$-module 
$A$ such that the following holds ?   
\[
A+f^*\left(\mathcal{E}_n\right)=\mathcal{R}_f.  
\]    
\end{question} 
\noindent 
Notice that by Ishikawa 
(\cite{ishikawainventiones, ishikawaaustralia}, 
see also \cite{ishikawa}), 
it is known if \lq\lq $f$ is finite and of corank one\rq\rq\; 
or \lq\lq it is $\mathcal{A}$-equivalent to 
a finite analytic map-germ\rq\rq, then 
there exists a finitely generated 
$f^*\left(\mathcal{E}_n\right)$-module $B$ 
satisfying the equality:  
\[
B+f^*\left(\mathcal{E}_n\right)=\mathcal{R}_f.  
\]
Notice also that in the case of Mather's $\mathcal{A}_e$ tangent space 
for a map-germ $g: (\mathbb{R}^n,0)\to (\mathbb{R}^p,0)$, 
the corresponding $\mathcal{E}_n$-module is 
nothing but $tg(\theta(n))$ (for the definition of Mather's 
$\mathcal{A}_e$ tangent space, see Section \ref{section 2}).    
Thus, Question \ref{question 1} asks whether or not 
the ramification module $\mathcal{R}_f$ has a similar structure 
as $T\mathcal{A}_e(g)$.        
\par 
\bigskip   
This paper is organized as  follows.   
In Section \ref{section 2}, 
preliminaries are given.   
Theorem \ref{theorem 1} (resp., Theorem \ref{theorem 2}) 
is proved in Section \ref{section 3} (resp., Section \ref{section 4}).   
Finally, in Section \ref{section 5}, examples concerning Theorem   
\ref{theorem 1} and Theorem \ref{theorem 2} are given.     
\section{Preliminaries}\label{section 2}
\subsection{Theory of singularities of mappings}\label{subsection 2.1}
\qquad 
In this subsection, it is partially 
reviewed several well-known notions/terminologies  
in the theory of singularities of mappings which are mainly developed by 
J.~Mather in \cite{mather3, mather4}.    
\cite{wall} is an excellent survey article on notions/terminologies reviewed 
in this subsection, which is recommended 
to readers.     
\par 
Two map-germ $f, g: (\mathbb{R}^n,0)\to \mathbb{R}^p$ 
are said to be $\mathcal{A}$-equivalent if there exist germs of 
diffeomorphism  
$h: (\mathbb{R}^n,0)\to (\mathbb{R}^n,0)$ and 
$H: (\mathbb{R}^p,f(0))\to (\mathbb{R}^p,g(0))$ such that 
$H\circ f\circ h=g$.       
\par 
Given a map-germ $f: (\mathbb{R}^n,0)\to (\mathbb{R}^p,0)$, 
the $\mathcal{A}$-equivalence class of $f$ is denoted by 
$\mathcal{A}(f)$.   
Let $\mathcal{E}_n$ be the $\mathbb{R}$-algebra consisting of 
function-germs $(\mathbb{R}^n,0)\to \mathbb{R}$ and 
let $m_n$ be the unique maximal ideal of $\mathcal{E}_n$.    
Let $f: (\mathbb{R}^n, 0)\to (\mathbb{R}^{p}, 0)$ be a 
map-germ.    
The set consisting of all vector fields along $f$ 
is denoted by $\theta(f)$.     
Notice that $\theta(f)$ is a finitely generated $\mathcal{E}_n$-module 
and it is a $\mathcal{E}_p$-module via $f$.    
For the identity map-germ  
$id_n: (\mathbb{R}^n, 0)\to (\mathbb{R}^{p}, 0)$, 
$\theta(id_n)$ is denoted by $\theta(n)$.  
The mapping $tf: \theta(n)\to \theta(f)$ 
(resp., $\omega f: \theta(p)\to \theta(f)$) is defined by 
$tf(\xi)=df\circ \xi$ (resp., $\omega f(\eta)=\eta\circ f$).    
The set $\mathcal{A}(f)$ may be regarded 
as an orbit of $f$ by the direct product 
of the following two groups:
\begin{eqnarray*}
\{\mbox{a germ of diffeomorphism } 
h: (\mathbb{R}^n,0)\to (\mathbb{R}^n,0)\},  \\ 
\{\mbox{a germ of diffeomorphism } 
H: (\mathbb{R}^p,0)\to (\mathbb{R}^p,0)\}.
\end{eqnarray*}   
Thus, the tangent space of $\mathcal{A}(f)$ at $f$ is naturally defined 
as follows:  
\[
T\mathcal{A}(f)=tf(m_n\theta(n))+\omega f(m_p\theta(p)).   
\]     
It is meaningful if $\{\mbox{diffeomorphism } 
h: (\mathbb{R}^n,0)\to (\mathbb{R}^n,0)\}$ 
(resp., $\{\mbox{diffeomorphism } H: (\mathbb{R}^p,0)\to 
(\mathbb{R}^p,0)\}$) is replaced with 
$\{\mbox{diffeomorphism } h: (\mathbb{R}^n,0)\to (\mathbb{R}^n,h(0))\}$ 
(resp., $\{\mbox{diffeomorphism } 
H: (\mathbb{R}^p,0)\to (\mathbb{R}^p,H(0))\}$).   
In this case, $tf(\theta(n))+\omega f(\theta(p))$ 
is denoted by $T\mathcal{A}_e(f)$.   
\[
T\mathcal{A}_e(f)=tf(\theta(n))+\omega f(\theta(p)). 
\]
\par 
Two map-germ $f, g: (\mathbb{R}^n,0)\to (\mathbb{R}^p,0)$ 
are said to be $\mathcal{R}$-equivalent if there exists a germ of 
diffeomorphism 
$h: (\mathbb{R}^n,0)\to (\mathbb{R}^n,0)$ such that 
$f\circ h=g$.   
Given a map-germ $f: (\mathbb{R}^n,0)\to (\mathbb{R}^p,0)$, 
similarly as in the case of $\mathcal{A}$-equivalence class of $f$ , 
$T\mathcal{R}(f)$ and  $T\mathcal{R}_e(f)$ can be naturally defined 
as follows:   
\begin{eqnarray*}
T\mathcal{R}(f) & = & tf(m_n\theta(n))  \\ 
T\mathcal{R}_e(f) & = & tf(\theta(n)).    
\end{eqnarray*}
Therefore, the condition (1) of the definition of frontal 
may be regarded as the condition 
that $\phi\in \theta(F)$ is perpendicular to $T\mathcal{R}_e(F)$ 
in $\theta(F)$.     
\par 
\smallskip 
Given a map-germ $f: (\mathbb{R}^n,0)\to (\mathbb{R}^p,0)$, 
the $\mathbb{R}$-algebra $\mathcal{E}_n/f^*m_p\mathcal{E}_n$ is denoted 
by $Q(f)$.     The $\mathbb{R}$-algebra $Q(f)$ is called 
the \textit{local algebra} of the map-germ $f$.    
Since the local algebra $Q(f)$ is an $\mathbb{R}$-algebra, it is a vector space.  
A map-germ $f: (\mathbb{R}^n,0)\to (\mathbb{R}^p,0)$ is said to be 
\textit{finite} if the vector space $Q(f)$ is of finite dimension.   
The dimension of $Q(f)$ as $\mathbb{R}$ vector space  is called 
the \textit{multiplicity} of $f$.    
\subsection{Openings}\label{subsection 2.2}
\qquad 
In this subsection, it is partially 
reviewed several well-known notions/terminologies  
in the theory of openings which are mainly developed by 
G.~Ishikawa.    
\cite{ishikawaaustralia, ishikawa} 
are excellent survey articles on notions/terminologies reviewed 
in this subsection, which are recommended 
to readers.     
There is one remark.   In \cite{ishikawa}, 
openings are defined for a map-germ 
$f: (\mathbb{R}^n, 0)\to (\mathbb{R}^m,0)$ $(n\le m)$; and 
in \cite{ishikawaaustralia} they are defined even for 
any multigerm $f: (\mathbb{R}^n, S)\to (\mathbb{R}^m, 0)$ $(n\le m)$.   
However, in this paper, 
it is needed only openings of an equidimensional 
mono-germ.    Thus, for the sake of clearness, in this subsection, 
we concentrate on reviewing notions/terminologies 
of opening only for an equidimensional mono-germ.   
\par 
\begin{definition}\label{versal opening}
{\rm 
Let $f: (\mathbb{R}^n,0)\to (\mathbb{R}^n,0)$ be 
an equidimensional map-germ.     
An opening $F=(f, \psi_1, \ldots, \psi_\ell)$ of $f$ is called a 
\textit{versal opening} (resp., \textit{mini-versal opening}) 
of $f$ if $1, \psi_1, \ldots, \psi_\ell$ form a system 
(resp., minimal system) of generators of $\mathcal{R}_f$ 
as an $f^*\left(\mathcal{E}_n\right)$-module.    
}
\end{definition}
\begin{lemma}[\cite{ishikawaaustralia}]
\label{lemma 2.1}
Let $f: (\mathbb{R}^n,0)\to (\mathbb{R}^n,0)$ be 
an equidimensional map-germ.   
Then, the following hold.   
\begin{enumerate}
\item[(1)]\quad 
$f^*\left(\mathcal{E}_n\right)
\subset \mathcal{R}_f\subset \mathcal{E}_n$.  
\item[(2)]\quad 
$\mathcal{R}_f$ is a $f^*\left(\mathcal{E}_n\right)$-module.   
\item[(3)]\quad 
$\mathcal{R}_f$ is a $C^\infty$ subring of $\mathcal{E}_n$.   
\item[(4)]\quad 
For a germ of diffeomorphism on the target space 
$H: (\mathbb{R}^n,0)\to (\mathbb{R}^n,0)$, 
$\mathcal{R}_{H\circ f}=\mathcal{R}_f$ holds.   
For a germ of diffeomorphism on the source space 
$h: (\mathbb{R}^n,0)\to (\mathbb{R}^n,0)$, 
$\mathcal{R}_{f\circ h}=h^*\left(\mathcal{R}_f\right)$ holds.  
\end{enumerate}
\end{lemma}
\begin{proposition}[\cite{ishikawainventiones}]
\label{corank one case}
Let $f: (\mathbb{R}^n,0)\to (\mathbb{R}^n,0)$ be 
a finite map-germ of corank one.   
Then, the following hold.  
\begin{enumerate}
\item[(1)]\quad 
$\mathcal{R}_f$ is a finitely generated $f^*\left(\mathcal{E}_n\right)$-module.   
Therefore, there is a versal opening of $f$.   
\item[(2)]\quad 
$1, \psi_1, \ldots, \psi_\ell\in \mathcal{R}_f$ generate $\mathcal{R}_f$ 
as $f(\mathcal{E}_n)$-module if and only if\;  
$[1], [\psi_1], \ldots, [\psi_\ell]$ generate the vector space 
$\mathcal{R}_f/\left(f^*m_n\right)\mathcal{R}_f$, 
where $[a]$ stands for $a+\left(f^*m_n\right)\mathcal{R}_f$.   
\end{enumerate}
\end{proposition}
\begin{proposition}
\label{analytic case}
Let $f: (\mathbb{R}^n,0)\to (\mathbb{R}^n,0)$ be 
a map-germ $\mathcal{A}$-equivalent to a finite 
analytic map-germ.   
Then, the following hold.  
\begin{enumerate}
\item[(1)]\quad 
$\mathcal{R}_f$ is a finitely generated $f^*\left(\mathcal{E}_n\right)$-module.   
Therefore, there is a versal opening of $f$.   
\item[(2)]\quad 
$1, \psi_1, \ldots, \psi_\ell\in \mathcal{R}_f$ generate $\mathcal{R}_f$ 
as $f(\mathcal{E}_n)$-module if and only if\;  
$[1], [\psi_1], \ldots, [\psi_\ell]$ generate the vector space 
$\mathcal{R}_f/\left(f^*m_n\right)\mathcal{R}_f$.    
\end{enumerate}
\end{proposition}
\noindent 
For the proof of Proposition \ref{analytic case}, see 
\cite{ishikawakodai, ishikawaaustralia}. 
\begin{proposition}[\cite{ishikawaaustralia}]
\label{injectivity}
Let $f: (\mathbb{R}^n,0)\to (\mathbb{R}^n,0)$ be a finite 
equidimensional map-germ.    
Then, every versal opening 
$F: (\mathbb{R}^n,0)\to (\mathbb{R}^{n+\ell},0)$ is injective.   
\end{proposition}
\noindent 
Here, a map-germ is said to be \textit{injective} if it has 
an injective representative.   
\section{Proof of Theorem \ref{theorem 1}}\label{section 3}
For any $i$ $(1\le i\le \ell)$, let the fundamental vector 
$(0, \ldots, 0, 1, 0, \ldots, 0)$ be denoted by ${\bf e}_i$.    
The vector ${\bf e}_i$ may be identified with the constant 
vector field $\frac{\partial}{\partial X_{n+i}}$ 
Let $\widetilde{Jf}$ be the cofactor matrix of 
the Jacobian matrix $Jf$.    
By using $\widetilde{Jf}$, 
for any $i$ $(1\le i\le \ell)$, 
define the map-germ $\Phi_i: (\mathbb{R}^n,0)\to T\mathbb{R}^{n+\ell}$ as 
$\Phi_i=(F, \phi_i)$, 
where 
\[
\phi_i=\left(
{ }^T\!\left(\left (|Jf|d\mu_i+2\mu_id|Jf|\right)\widetilde{Jf}\right), 
-{\bf e}_i
\right).   
\]  
Here, $d\varphi$ for a function-germ 
$h: (\mathbb{R}^n,0)\to \mathbb{R}$ (resp., ${ }^T{\bf v}$ for a 
vector field ${\bf v}$) stands for  the gradient vector field of $h$ 
(resp., the transposed vector field of  
${\bf v}$).   
Then, it follows that $\Phi_j$ is of class $C^\infty$.   
Thus, $\Phi_j$ is a $C^\infty$ vector field along $F$.   
It is clear that $\Phi_1, \ldots, \Phi_\ell$ satisfy the condition 
(2) and (3) of the definition of frontal.   
Therefore, it is sufficient to show that $\Phi_i$ satisfies the condition 
(1) for any $i$ $(1\le i\le \ell)$.     
Let $\xi$ be an element of $\theta(n)$.     
Then, 
\begin{eqnarray*}
{ } & { } & \phi_i(x)\cdot tF(\xi)(x) \\ 
{ } & = &  \left(|Jf|(x)d\mu_i(x)+2\mu_i(x)d|Jf|(x)\right) 
\widetilde{Jf}(x)Jf(x)\xi(x) \\ 
{ } & { } & \qquad\qquad 
- 
|Jf|(x){ } \left(|Jf|(x)d\mu_i(x)+2\mu_i(x)d|Jf|(x)\right)\xi(x) \\ 
{ } & = & 
|Jf|(x)\left(\left(|Jf|(x)d\mu_i(x)+2\mu_i(x)d|Jf|(x)\right)\xi(x)\right.\\ 
{ } & { } & \qquad\qquad 
-\left.\left(|Jf|(x)d\mu_i(x)+2\mu_i(x)d|Jf|(x)\right)\xi(x)\right) \\ 
{ } & = & 0.   
\end{eqnarray*}
Hence, $F$ is a frontal.   
\hfill $\Box$ 
\section{Proof of Theorem \ref{theorem 2}}
\label{section 4}
\subsection{The case $n=1$}\label{subsection 4.1}
\quad
Set $\delta=\dim_{\mathbb{R}}Q(f)$.    
In the case of $n=1$, there exists a germ of diffeomorphism 
$h: (\mathbb{R}, 0)\to (\mathbb{R},0)$ such that 
$g(x)=f\circ h(x)=\frac{1}{\delta}x^\delta$.     
By \cite{ishikawatransaction}, 
the ramification module of $g$ can be described as follows.   
\[
\mathcal{R}_g = 
\left\{
\begin{array}{ll}
\left\langle 1 \right\rangle_{g^*\left(\mathcal{E}_1\right)}
=\left\langle 1 \right\rangle_{\mathcal{E}_1}\quad 
& (\delta=1), \\ 
\left\langle 1, x^3 \right\rangle_{g^*\left(\mathcal{E}_1\right)}
\quad 
& (\delta=2), \\ 
\left\langle 1, x^4, x^5 \right\rangle_{g^*\left(\mathcal{E}_1\right)}
\quad 
& (\delta=3). \\  
\end{array}
\right.
\]
On the other hand, since $|Jg|^2(x)=x^{2(\delta-1)}$, 
the following holds. 
\[
\left\langle |Jg|^2 \right\rangle_{\mathcal{E}_1} = 
\left\{
\begin{array}{ll}
\left\langle 1 \right\rangle_{\mathcal{E}_1}\quad 
& (\delta=1), \\ 
\left\langle x^2, x^3 \right\rangle_{\mathcal{E}_1}
\quad 
& (\delta=2), \\ 
\left\langle x^4, x^5 \right\rangle_{\mathcal{E}_1}
\quad 
& (\delta=3). \\  
\end{array}
\right.
\]   
Thus, we have the following inclusion.  
\[
\left\langle |Jg|^2 \right\rangle_{\mathcal{E}_1}+ 
g^*\left(\mathcal{E}_1\right) 
\supset \mathcal{R}_g.  
\]
Combining this inclusion with Theorem \ref{theorem 1}, 
the following equality holds.   
\[
\left\langle |Jg|^2 \right\rangle_{\mathcal{E}_1}+ 
g^*\left(\mathcal{E}_1\right) 
= \mathcal{R}_g. 
\]
Since $h$ is a germ of diffeomorphism, we have the following by the chain rule.  
\[
\left\langle |Jg|^2 \right\rangle_{\mathcal{E}_1}
= 
\left\langle |J(f\circ h)|^2 \right\rangle_{\mathcal{E}_1}
= 
h^*\left(\left\langle |Jf|^2 \right\rangle_{\mathcal{E}_1}\right).   
\]
It follows 
\[
h^*\left(\left\langle |Jf|^2 \right\rangle_{\mathcal{E}_1}+ 
f^*\left(\mathcal{E}_1\right)\right) 
= h^*\left(\mathcal{R}_f\right). 
\]
Hence, by pulling both sides back by $\left(h^{-1}\right)^*$, 
the desired equality can be obtained as follows.   
\[
\left\langle |Jf|^2 \right\rangle_{\mathcal{E}_1}+ 
f^*\left(\mathcal{E}_1\right) 
= \mathcal{R}_f. 
\]
\hfill 
$\Box$
\subsection{The case $n\ge 2$}\label{subsection 4.2}
\quad 
Again in this subsection, $\dim_{\mathbb{R}}Q(f)$ is denoted by $\delta$.   
The assumption $\delta\le 3$ implies that $f$ is of corank one.   
Since $\frac{1}{2}x^2+t x$ (resp., $\frac{1}{3}x^3+t x$) is 
an $\mathcal{R}$-versal unfolding of 
$\frac{1}{2}x^2$ (resp., $\frac{1}{3}x^3$), it is deduced that 
there exist germs of diffeomorphism 
$h, H: (\mathbb{R}^n,0)\to (\mathbb{R}^n,0)$ and 
a function-germ $\alpha: (\mathbb{R}^{n-1}, 0)\to \mathbb{R}$ 
such that 
\[
H\circ f\circ h (x, \lambda_1, \ldots, \lambda_{n-1}) 
= 
\left\{
\begin{array}{ll}
\left(x, \lambda_1, \ldots, \lambda_{n-1}\right)\quad 
& (\delta=1), \\ 
\left(\frac{1}{2}x^2+\alpha(\lambda_1, \ldots, \lambda_{n-1})x, 
\lambda_1, \ldots, \lambda_{n-1}\right)\quad 
& (\delta=2), \\ 
\left(\frac{1}{3}x^3+\alpha(\lambda_1, \ldots, \lambda_{n-1})x, 
\lambda_1, \ldots, \lambda_{n-1}\right)\quad 
& (\delta=3),    
\end{array}
\right.
\]
where $(x, \lambda_1, \ldots, \lambda_{n-1})$ is an element of 
$\mathbb{R}\times \mathbb{R}^{n-1}=\mathbb{R}^n$.    
Set $g=H\circ f\circ h$ and 
$\lambda=\left(\lambda_1, \ldots, \lambda_{n-1}\right)$.       
By \cite{ishikawatransaction}, 
the ramification module of $g$ can be described as follows.   
\[
\mathcal{R}_g = 
\left\{
\begin{array}{ll}
\left\langle 1 \right\rangle_{g^*\left(\mathcal{E}_n\right)}
=\left\langle 1 \right\rangle_{\mathcal{E}_n}\quad 
& (\delta=1), \\ 
\left\langle 1, \frac{1}{3}x^3
+\frac{1}{2}\alpha(\lambda)x^2  
\right\rangle_{g^*\left(\mathcal{E}_n\right)}
\quad 
& (\delta=2), \\ 
\left\langle 1, \frac{1}{4}x^4+\frac{1}{2}\alpha(\lambda)x^2, 
\frac{1}{5}x^5+\frac{1}{3}\alpha(\lambda)x^3 \right\rangle_{g^*\left(\mathcal{E}_n\right)}
\quad 
& (\delta=3). \\  
\end{array}
\right.
\]
On the other hand, we have  
\[
\left\langle |Jg|^2\right\rangle_{\mathcal{E}_n} = 
\left\{
\begin{array}{ll}
\left\langle 1 \right\rangle_{\mathcal{E}_n}\quad 
& (\delta=1), \\ 
\left\langle (x+\alpha(\lambda))^2 \right\rangle_{\mathcal{E}_n}
\quad 
& (\delta=2), \\ 
\left\langle (x^2+\alpha(\lambda))^2 \right\rangle_{\mathcal{E}_n}
\quad 
& (\delta=3). \\  
\end{array}
\right.
\]   
Since 
\begin{eqnarray*}
\left(x+\alpha(\lambda)\right)^2 
& = &  
2\left(\frac{1}{2}x^2+\alpha(\lambda)x\right)+ 
\left(\alpha(\lambda)\right)^2, \\ 
x\left(x+\alpha(\lambda)\right)^2 
& = &  
3\left(\frac{1}{3}x^3+\frac{1}{2}\alpha(\lambda)x^2\right)+ 
\alpha(\lambda)\left(\frac{1}{2}x^2+\alpha(\lambda)x\right), \\ 
\left(x^2+\alpha(\lambda)\right)^2 
& = &  
4\left(\frac{1}{4}x^4+\frac{1}{2}\alpha(\lambda)x^2\right)+ 
\left(\alpha(\lambda)\right)^2\quad \mbox{and} \\  
x\left(x^2+\alpha(\lambda)\right)^2 
& = &  
5\left(\frac{1}{5}x^5+\frac{1}{3}\alpha(\lambda)x^3\right)+ 
\alpha(\lambda)\left(\frac{1}{3}x^3+\alpha(\lambda)x\right), 
\end{eqnarray*}
the following holds.   
\[
\left\langle |Jg|^2 \right\rangle_{\mathcal{E}_n}+ 
g^*\left(\mathcal{E}_n\right) 
\supset \mathcal{R}_g.  
\]
Combining this inclusion with Theorem \ref{theorem 1}, 
we have the following equality.   
\[
\left\langle |Jg|^2 \right\rangle_{\mathcal{E}_n}+ 
g^*\left(\mathcal{E}_n\right) 
= \mathcal{R}_g. 
\]
Since $h, H$ are germs of diffeomorphism, it follows 
the following.  
\[
\left\langle |Jg|^2 \right\rangle_{\mathcal{E}_n}
= 
\left\langle |J(H\circ f\circ h)|^2 \right\rangle_{\mathcal{E}_n}
= 
h^*\left(\left\langle |Jf|^2 \right\rangle_{\mathcal{E}_n}\right).   
\]
Combining this property with two facts  
$g^*\left(\mathcal{E}_n\right)=
h^*\left((H\circ f)^*\left(\mathcal{E}_n\right)\right)
=h^*\left(f^*\left(\mathcal{E}_n\right)\right)
$ 
and 
$\mathcal{R}_g=\mathcal{R}_{H\circ f\circ h}= 
h^*\left(\mathcal{R}_{H\circ f}\right)
=h^*\left(\mathcal{R}_f\right)
$, 
we have  
\[
h^*\left(\langle |Jf|^2 \rangle_{\mathcal{E}_n}+ 
f^*\left(\mathcal{E}_n\right)\right) 
= h^*\left(\mathcal{R}_f\right). 
\]
Therefore, also in this case, 
we have the following desired equality.   
\[
\left\langle |Jf|^2 \right\rangle_{\mathcal{E}_n}+ 
f^*\left(\mathcal{E}_n\right) 
= \mathcal{R}_f. 
\]
\hfill 
$\Box$
\section{Examples of Theorems \ref{theorem 1} and 
\ref{theorem 2}}\label{section 5}
All map-germs $h_i: (\mathbb{R}^n,0)\to (\mathbb{R}^n,0)$ and 
$H_j: (\mathbb{R}^{n+\ell},0)\to (\mathbb{R}^{n+\ell},0)$ 
appearing in this section are 
germs of diffeomorphism.     
\subsection{Fold Singularity} \label{fold singularity}
Let $f: (\mathbb{R}^2, 0)\to (\mathbb{R}^2,0)$ be 
the map-germ defined by 
$f(x,y)=\left(\frac{1}{2}x^2+xy,\; y\right)$. 
The map-germ $f$ is $\mathcal{A}$-equivalent to the map-germ 
named \textit{Type 2} in the list of \cite{rieger}.    
It is clear that $|Jf|(x,y)=x+y$ and $\dim_{\mathbb{R}}Q(f)=2$.   
Set $\mu_1(x,y)=1$.     
Then, by Theorem \ref{theorem 1}, 
the map-germ $F: (\mathbb{R}^2,0)\to (\mathbb{R}^3,0)$ 
defined by 
\begin{eqnarray*}
F(x,y)  & =  & \left(f(x,y), \mu_1(x,y)|Jf|^2(x,y)\right) \\ 
{ } & = &  \left(\frac{1}{2}x^2+xy,\; y,\; (x+y)^2\right) =  
\left(\frac{1}{2}\left(x+y\right)^2-\frac{1}{2}y^2,\; y,\; (x+y)^2\right)    
\end{eqnarray*}
is a frontal.   
Set $H_1(X,Y,Z)=(X-\frac{1}{2}Y^2,\; Y,\; Z)$.    
Then, 
\[
H_1\circ F(x,y)  = 
\left(\frac{1}{2}(x+y)^2, \; y, \; 
(x+y)^2\right).    
\]   
Secondly, set $h_1(x,y)=\left(x-y, y\right)$ and 
$H_2(X,Y,Z)=(2X,\; Y,\; Z-2X)$.   Then, 
\[
H_2\circ H_1\circ F\circ h_1(x,y)  
= 
\left(x^2, \; y, \; 0\right).    
\]
In Differential Geometry, 
the map-germ $H_2\circ H_1\circ F\circ h_1$ is called 
the normal form of {\it fold singularity} (\cite{fujimorietal}). 
\subsection{Cuspidal Edge}\label{cuspidal edge}
Not only in this subsection, but also 
until Subsection \ref{open folded umbrella}, 
we start from the same map-germ 
$f: (\mathbb{R}^2, 0)\to (\mathbb{R}^2,0)$ as in 
Subsection \ref{fold singularity}.   
Namely, in this subsection, $f$ is the map-germ defined by 
$
f(x,y)=\left(\frac{1}{2}x^2+xy,\; y\right) 
$. 
Set $\mu_1(x,y)=x$.     
Then, by Theorem \ref{theorem 1}, 
the map-germ $F: (\mathbb{R}^2,0)\to (\mathbb{R}^3,0)$ 
defined by 
\begin{eqnarray*}
F(x,y) & = & \left(f(x,y), \mu_1(x,y)|Jf|^2(x,y)\right)\\ 
{ } & = & \left(\frac{1}{2}x^2+xy,\; y,\; x(x+y)^2\right) \\ 
{ } & = &  \left(\frac{1}{2}\left(x+y\right)^2-\frac{1}{2}y^2,
\; y,\; (x+y)^3-y(x+y)^2\right)    
\end{eqnarray*}
is a frontal.   
Set $H_1(X,Y,Z)=(X-\frac{1}{2}Y^2,\; Y,\; Z)$.    
Then, 
\begin{eqnarray*}
H_1\circ F(x,y) 
& = &  \left(\frac{1}{2}(x+y)^2, \; y, \; 
(x+y)^3-y(x+y)^2\right).    
\end{eqnarray*}
Secondly, set $h_1(x,y)=\left(x-y, y\right)$ and 
$H_2(X,Y,Z)=(2X,\; Y,\; Z-2XY)$.   Then, 
\[
H_2\circ H_1\circ F\circ h_1(x,y)  
= 
\left(x^2, \; y, \; x^3\right),     
\]
well-known as the normal form of 
\textit{cuspidal edge} (for instance, see 
\cite{krsuy}). 
\subsection{Folded Umbrella(Cuspidal Crosscap)} \label{folded umbrella}
As explained in the last subsection, in this subsection again, 
$f: (\mathbb{R}^2, 0)\to (\mathbb{R}^2,0)$ is  
the map-germ defined by $f(x,y)=\left(\frac{1}{2}x^2+xy,\; y\right)$. 
Set $\mu_1(x,y)=x^2$.     
Then, by Theorem \ref{theorem 1}, 
the map-germ $F: (\mathbb{R}^2,0)\to (\mathbb{R}^3,0)$ 
defined by 
\begin{eqnarray*}
F(x,y) & = & \left(f(x,y), \mu_1(x,y)|Jf|^2(x,y)\right)\\ 
{ } & = & \left(\frac{1}{2}x^2+xy,\; y,\; x^2(x+y)^2\right) \\ 
{ } & = &  \left(\frac{1}{2}x^2+xy,\; y,\; x^4+2x^3y+x^2y^2\right)    
\end{eqnarray*}
is a frontal.   
Set $H_1(X,Y,Z)=(X,\; Y,\; Z-X^2)$.    
Then, 
\begin{eqnarray*}
H_1\circ F(x,y) 
& = &  \left(\frac{1}{2}x^2+xy, \; y, \; 
\frac{3}{4}x^4+x^3y\right).    
\end{eqnarray*}
Secondly, set $h_1(x,y)=\left(x, \frac{1}{2}y\right)$.   Then, 
\[
H_1\circ F\circ h_1(x,y)  
= 
\left(\frac{1}{2}x^2+\frac{1}{2}xy, \; \frac{1}{2}y, \; 
\frac{3}{4}x^4+\frac{1}{2}x^3y\right).    
\]
Thirdly, set 
$H_2(X, Y, Z)=\left(2X, 2Y, \frac{4}{3}Z\right)$.   Then, 
we have 
\[
H_2\circ H_1\circ F\circ h_1(x,y)  
= 
\left(x^2+xy, \; y, \; 
x^4+\frac{2}{3}x^3y\right),     
\]
called the normal form of 
{\it folded umbrella} (for instance, see \cite{ishikawa}). 
A folded umbrella is also called a 
{\it cuspidal crosscap}.   
As the normal form of cuspidal crosscap, 
Some references adopt 
$(x,y)\mapsto \left(x^2,\; y, \; x^3y\right)$ 
(for instance, see \cite{fujimorisajiumeharayamada}).    
Set $h_2(x,y)=\left(\frac{1}{2}(x-y), y\right)$ 
and 
\begin{eqnarray*}
{ } & { } & 
H_3(X,Y,Z)\\ 
{ } & = & \left(4\left(X+\frac{1}{4}Y^2\right), \; Y, \; 
-6\left(Z-\left(X+\frac{1}{4}Y^2\right)^2
-\frac{1}{2}\left(X+\frac{1}{4}Y^2\right)Y^2
- \frac{1}{16}Y^4\right)\right).   
\end{eqnarray*}   
Then, 
it is easily confirmed 
\[
H_3\circ H_2\circ H_1\circ F\circ h_1\circ h_2(x,y)  
 = 
\left(x^2,\; y,\; x^3y\right).   
\]
\subsection{Open Folded Umbrella}\label{open folded umbrella}
As explained, we start from 
the same $f: (\mathbb{R}^2, 0)\to (\mathbb{R}^2,0)$ 
as in Subsection \ref{fold singularity}.   
Thus, $f$ is defined by $f(x,y)=\left(\frac{1}{2}x^2+xy,\; y\right)$ 
in this subsection. 
Suppose that $\ell \ge 3$ and 
set $\mu_1(x,y)=x^2, \mu_2(x,y)=x^2|Jf|(x,y)$ and 
$\mu_i(x,y)=0$ $(3\le i\le \ell)$.    
Then, 
\begin{eqnarray*}
F(x,y) & = & \left(f(x,y), \mu_1|Jf|^2(x,y), \ldots, 
\mu_\ell(x,y)|Jf|^2(x,y)\right) \\  
{ }  & = & 
\left(\frac{1}{2}x^2+xy,\; y,\; x^2(x+y)^2, \; x^2(x+y)^3, 0, \ldots, 0\right) \\ 
{ } & = &  \left(\frac{1}{2}x^2+xy,\; y,\; x^4+2x^3y+x^2y^2, \right. \\ 
{ } & { } & \qquad \left. x^5+3x^4y+3x^3y^2+x^2y^3, 
0, \ldots, 0\right).     
\end{eqnarray*}
Set 
\[
H_1(X,Y,U_1, U_2, U_3, \ldots, U_\ell)
=(X,\; Y,\; U_1-X^2,\; U_2-YU_1,\; U_3, \ldots, U_\ell).
\]    
Then, 
\begin{eqnarray*}
H_1\circ F(x,y) 
& = &  \left(\frac{1}{2}x^2+xy, \; y, \; 
\frac{3}{4}x^4+x^3y, \; x^5+2x^4y+x^3y^2, 0, \ldots, 0\right).    
\end{eqnarray*}
Set $h_1(x,y)=\left(x, \frac{1}{2}y\right)$.   Then, 
\begin{eqnarray*}
H_1\circ F\circ h_1(x,y)  
& = &  
\left(\frac{1}{2}x^2+\frac{1}{2}xy, \; \frac{1}{2}y, \; 
\frac{3}{4}x^4+\frac{1}{2}x^3y, \; 
x^5+x^4y+\frac{1}{4}x^3y^2, 0, \ldots, 0\right).      
\end{eqnarray*}
Nextly, set 
\[
H_2(X, Y, U_1, U_2, U_3, \ldots, U_\ell)=
\left(X,\; Y,\; U_1,\; U_2-YU_1,\; U_3, \ldots, U_\ell\right).
\]   Then, 
we have 
\begin{eqnarray*}
{ } & { } & H_2\circ H_1\circ F\circ h_1(x,y) \\   
{ } & = &  
\left(\frac{1}{2}x^2+\frac{1}{2}xy, \; \frac{1}{2}y, \; 
\frac{3}{4}x^4+\frac{1}{2}x^3y, \; x^5+\frac{5}{8}x^4y, 
0, \ldots, 0\right).       
\end{eqnarray*}
Finally, set 
\[
H_3(X, Y, U_1, U_2, \ldots, U_\ell)= 
\left(2X,\; 2Y,\; \frac{4}{3}U_1,\; U_2, \ldots, U_\ell\right).
\]
Then, 
\begin{eqnarray*}
{ } & { } & H_3\circ H_2\circ H_1\circ F\circ h_1(x,y) \\ 
{ } & = & 
\left(x^2+xy, \; y, \; 
x^4+\frac{2}{3}x^3y,\; x^5+\frac{5}{8}x^4y,\;  
0, \ldots, 0
\right), 
\end{eqnarray*}
called the normal form of {\it open folded umbrella} (\cite{arnold2, ishikawa}).   
\subsection{Swallowtail}\label{swallowtail}
Let $f: (\mathbb{R}^2, 0)\to (\mathbb{R}^2,0)$ be 
the map-germ defined by $f(x,y)=\left(\frac{1}{3}x^3+xy,\; y\right)$. 
The map-germ $f$ is $\mathcal{A}$-equivalent to the map-germ 
named \textit{Type 3} in the list of \cite{rieger}.    
It is clear that $|Jf|(x,y)=x^2+y$ and $\dim_{\mathbb{R}}Q(f)=3$.   
Set $\mu_1(x,y)=1$.    
We consider the map-germ $F: (\mathbb{R}^2,0)\to (\mathbb{R}^3,0)$ 
defined by 
\begin{eqnarray*}
F(x,y) & = & \left(f(x,y), \mu_1(x,y)|Jf|^2(x,y)\right)\\ 
{ } & = & \left(\frac{1}{3}x^3+xy,\; y,\; (x^2+y)^2\right) \\ 
{ } & = &  \left(\frac{1}{3}x^3+xy,\; y,\; x^4+2x^2y+y^2\right).    
\end{eqnarray*}
By Theorem \ref{theorem 1}, $F$ is a frontal.  
Set $H_1(X,\; Y,\; Z)=(X,\; Y,\; Z-Y^2)$, 
$H_2(X,\;Y,\; Z)=(-12X,\; 6Y,\; 3Z)$.   
Then, 
\[
H_2\circ H_1\circ {F}(x, y)=
\left(-4x^3-12xy,6y, 3x^4+6x^2y \right). 
\]
Next, set $h_1(x,\; y)=\left(x,\; \frac{1}{6}y\right)$.     Then, 
\[
H_2\circ H_1\circ \widetilde{F}\circ h_1(x,\; y)= 
\left(-4x^3-2xy,\; y, \; 3x^4+x^2y\right), 
\]
well-known as the normal form of 
{\it swallowtail} (for instance see \cite{brucegiblin}, 129 page).  
\subsection{Open Swallowtail}\label{open swallowtail}
As in Subsection \ref{swallowtail}, 
let $f: (\mathbb{R}^2, 0)\to (\mathbb{R}^2,0)$ 
be the map-germ defined by $f(x,y)=\left(\frac{1}{3}x^3+xy,\; y\right)$. 
Suppose that $\ell\ge 3$ and 
set $\mu_1(x,y)=1, \mu_2(x,y)=x, \mu_i(x,y)=0$ $(3\le i\le \ell)$.     
Let $F: (\mathbb{R}^2, 0)\to (\mathbb{R}^{2+\ell},0)$ be the map-germ 
defined by  
\begin{eqnarray*}
F(x,y) & = & \left(f(x,y), \mu_1|Jf|^2(x,y), \ldots, \mu_\ell|Jf|^2(x,y)\right) \\ 
{ } & = & \left(\frac{1}{3}x^3+xy,\; y,\; (x^2+y)^2,\; x(x^2+y)^2,\; 
0, \ldots, 0\right) \\ 
{ } & = &  \left(\frac{1}{3}x^3+xy,\; y,\; x^4+2x^2y+y^2, 
x^5+2x^3y+xy^2, 0, \ldots, 0\right).    
\end{eqnarray*}
Set 
\[
H_1(X,Y,U_1, U_2, U_3, \ldots, U_\ell)
=(X,\; Y,\; U_1-Y^2,\; U_2-XY,\; U_3, \ldots, U_\ell).  
\]
Then, 
\[
H_1\circ F(x,y) 
 =   \left(\frac{1}{3}x^3+xy, \; y, \; 
x^4+2x^2y,\;x^5+\frac{5}{3}x^3y,\; 0, \ldots, 0\right).    
\]
Set $h_1(x,y)=\left(x,\; \frac{1}{3}y\right)$.   Then, 
\[ 
H_1\circ F\circ h_1(x,y)  
 =   \left(\frac{1}{3}x^3+\frac{1}{3}xy, \; \frac{1}{3}y, \; 
x^4+\frac{2}{3}x^2y\; x^5+\frac{5}{9}x^3y,\; 0, \ldots, 0\right).    
\] 
Set $H_2(X,Y,U_1, \ldots, U_\ell)
=(3X\; ,3Y,\; U_1, \ldots, U_\ell)   
$.   Then, 
\[
H_2\circ H_1\circ G\circ h_1(x,y) 
 =   \left(x^3+xy, \; y, \; 
x^4+\frac{2}{3}x^2y, x^5+\frac{5}{9}x^3y,\; 0, 
\ldots, 0\right), 
\] 
known as the normal form of 
{\it open swallowtail} (\cite{ishikawa}).   
\subsection{Opening of $4_k$ in the list of \cite{rieger}}
\label{opening of $4_k$}  
Let $k$ be an integer greater than $1$ and let 
$f_{k, \pm}: (\mathbb{R}^2, 0)\to (\mathbb{R}^2,0)$ be 
the map-germ defined by 
$f_{k, \pm}(x,y)=\left(\frac{1}{3}x^3\pm xy^k,\; y\right)$. 
The map-germ $f_{k, \pm}$ is $\mathcal{A}$-equivalent to the map-germ 
named \textit{Type $4_k$} in the list of \cite{rieger}.    
It is clear that $|Jf|(x,y)=x^2\pm y^k$ and $\dim_{\mathbb{R}}Q(f)=3$.   
Set $\mu_1(x,y)=1$.    
We consider the map-germ 
$F_{k, \pm}: (\mathbb{R}^2,0)\to (\mathbb{R}^3,0)$ 
defined by 
\begin{eqnarray*}
F_{k, \pm}(x,y) & = & \left(f_{k, \pm}(x,y), \mu_1(x,y)|Jf|^2(x,y)\right)\\ 
{ } & = & \left(\frac{1}{3}x^3\pm xy^k,\; y,\; (x^2\pm y^k)^2\right) \\ 
{ } & = &  \left(\frac{1}{3}x^3\pm xy^k,\; y,\; x^4\pm 2x^2y^k+y^{2k}\right).    
\end{eqnarray*}
Theorem \ref{theorem 1} guarantees that 
$F_{k, \pm}$ is a frontal.  
Set $H_1(X,\; Y,\; Z)=(X,\; Y,\; Z-Y^{2k})$ and  
$H_2(X,\;Y,\; Z)=(6X,\; 6Y,\; 3Z)$.   
Then, 
\[
H_2\circ H_1\circ {F_{k, \pm}}(x, y)=
\left(2x^3\pm 6xy^k,\; 6y,\; 3x^4\pm 6x^2y^k\right). 
\]
Next, set $h_1(x,\; y)=\left(x,\; \frac{1}{\sqrt[k]{6}}y\right)$  
and $H_3(X,\; Y,\; Z)=(X,\; \frac{\sqrt[k]{6}}{6}Y,\; Z)$.     
Then, 
\[
H_3\circ H_2\circ H_1\circ {F_{k,\pm}}\circ h_1(x,\; y)= 
\left(2x^3\pm xy^k,\; y,\; 3x^4\pm x^2y^k\right).   
\]
The form of  
$H_3\circ H_2\circ H_1\circ {F_{2, -}}\circ h_1$ is 
exactly the same as the map-germ called a 
\textit{double swallowtail} given in \cite{sajiumeharayamada}.  
\subsection{A non-analytic example}
\label{non-analytic example}  
Let $\psi: (\mathbb{R},0)\to (\mathbb{R},0)$ be a flat function-germ 
(i.e. $j^\infty\psi(0)=0$) and let 
$f_{\infty}: (\mathbb{R}^2, 0)\to (\mathbb{R}^2,0)$ 
be the map-germ defined by 
$f_{\infty}(x,y)=\left(\frac{1}{3}x^3\pm x\psi(y),\; y\right)$. 
It is clear that $\dim_{\mathbb{R}}Q(f_\infty)=3$.   
Thus, by Theorem \ref{theorem 2}, 
the following equality holds even for the $f_\infty$.    
\[
\left\langle |Jf_\infty|^2 \right\rangle_{\mathcal{E}_2}+
f^*_{\infty}\left(\mathcal{E}_2\right)=\mathcal{R}_{f_\infty}.    
\]
Therefore, for any opening 
$F: (\mathbb{R}^2, 0)\to (\mathbb{R}^{2+\ell},0)$ of $f_\infty$ 
there exist germs of diffeomorphisms 
$h: (\mathbb{R}^2, 0)\to (\mathbb{R}^2,0)$, 
$H: (\mathbb{R}^{2+\ell}, 0)\to (\mathbb{R}^{2+\ell},0)$ 
and function-germs 
$\mu_1, \ldots, \mu_\ell: (\mathbb{R}^2,0)\to \mathbb{R}$ 
such that the following holds.   
\[
H\circ F\circ h=
(f_\infty,\; \mu_1|Jf_\infty|^2,\; \ldots,\; \mu_\ell|Jf_\infty|^2).   
\]
\subsection{$A_k$-Front Singularity}\label{$A_k$-front singularity}
Let $k$ be an integer greater than $1$ and let 
$f_k: (\mathbb{R}^k,0)\to (\mathbb{R}^k,0)$ be the map-germ 
defined as follows:
\[
f_k(x_1, \ldots, x_k)  
 =  
\left(\frac{1}{k+1}x_1^{k+1}+\frac{1}{k-1}x_1^{k-1}x_2
+\cdots + x_1x_k, x_2, \ldots, x_k\right).   
\]
The map-germ $f_k$ is a generalization of $f$ 
given in Subsection \ref{swallowtail}  
because $f_2$ is exactly the same as the $f$.        
The map-germ $f_k$ is well-known as the normal form of 
corank-one isolated stable singularity 
(for instance, see \cite{arnoldguseinzadevarchenko, mather4}).  
For the $f_k$,  
$\dim_{\mathbb{R}}Q(f_k)=k+1$ and 
$|Jf_k|^2(x_1,0,\ldots, 0)=x_1^k$.     
\par 
From now on, suppose that $k$ is greater than $2$.     
Let $\mu: (\mathbb{R}^n, 0)\to \mathbb{R}$ be a function-germ.   
Since $2k-(k+1)=k-1>1$, for the normal vector field 
$\nu$ of the frontal 
$(f_k, \mu|Jf_k|^2)$, $\nu(x_1, 0, \ldots, 0)$ must be singular at the origin.   
\par 
On the other hand, since $f_k$ is a polynomial map-germ of corank one, 
$\mathcal{R}_{f_k}$ can be described explicitely 
by \cite{ishikawatransaction}.     
In particular, there exists an opening $\widetilde{F}_k$ of $f_k$ whose 
normal vector field $\nu$ is non-singular with respect to 
$x_1$ (A concrete construction of $\nu$ can be found 
in \cite{sajiumeharayamadacambridge}).     
Therefore, if $k\ge 3$, it is impossible to expect 
the converse of Theorem \ref{theorem 1} 
for the frontal $\widetilde{F}_k$.       
\section*{Acknowledgement}
The author would like to thank Goo Ishikawa 
for his helpful suggestions.   
This work was 
supported 
by JSPS KAKENHI Grant Number 17K05245.   

\end{document}